\documentclass[a4paper, final, 11pt]{article}

\usepackage{amsfonts,amsbsy,enumerate} 
\usepackage{enumerate}
\usepackage{graphicx,psfrag}
\newtheorem{proposition}{Proposition}[section]
\newtheorem{theorem}[proposition]{Theorem}
\newtheorem{lemma}[proposition]{Lemma}

\newenvironment{proof}{\smallskip\noindent\emph{\textbf{Proof.}}\hspace{1pt}}%
{\hspace{-5pt}{\nobreak\quad\nobreak\hfill\nobreak$\square$\vspace{8pt}%
\par}\smallskip\goodbreak}

\newenvironment{proofof}[1]{\smallskip\noindent\emph{\textbf{Proof of #1.}}%
\hspace{1pt}}{\hspace{-5pt}{\nobreak\quad\nobreak\hfill\nobreak%
$\square$\vspace{8pt}\par}\smallskip\goodbreak}

\newcommand{\Section}[1]{\section{#1}\setcounter{equation}{0}}

\newcommand{\Id}{\mathinner{\mathrm{Id}}}

\newcommand{\comp}{\mathop\bigcirc}

\newcommand{\C}[1]{\mathbf{C^{#1}}}

\newcommand{\modulo}[1]{{\left|#1\right|}}
\newcommand{\norma}[1]{{\left\|#1\right\|}}
\newcommand{\Ref}[1]{{\rm(\ref{#1})}}
\newcommand{\reali}{{\mathbb{R}}}

\newcommand{\naturali}{{\mathbb{N}}}
\newcommand{\tv}{\mathrm{TV}}

\renewcommand{\O}{\mathinner{\mathcal{O}(1)}}
\renewcommand{\epsilon}{\varepsilon}
\renewcommand{\phi}{\varphi}
\renewcommand{\theta}{\vartheta}
\renewcommand{\L}[1]{\mathbf{L^#1}}

\newcommand{\caratt}[1]{{\chi_{\strut#1}}}
\newcommand{\gf}{{\mathinner{\mathbf{\Phi}}}}

\newcommand{\F}{{\hat F}}

\title{Hyperbolic Balance Laws \\ with a \\ Dissipative Non Local Source}

\author{Rinaldo M.~Colombo \\ \small Dipartimento di Matematica \\
  \small Universit\`a degli Studi di Brescia \\ \small Via Branze, 38
  \\ \small 25123 Brescia, Italy \\ \texttt{Rinaldo.Colombo@UniBs.it}\\
  \and Graziano Guerra \\ \small Dip.~di Matematica e Applicazioni \\
  \small Universit\`a di Milano -- Bicocca \\ \small Via Bicocca degli
  Arcimboldi, 8 \\ \small 20126 Milano, Italy \\
  \texttt{Graziano.Guerra@UniMiB.it}}

\begin{document}

\maketitle

\begin{abstract}

  \noindent This paper considers systems of balance law with a
  dissipative non local source. A global in time well posedness result
  is obtained. Estimates on the dependence of solutions from the flow
  and from the source term are also provided. The technique relies on
  a recent result on quasidifferential equations in metric spaces.

  \medskip

  \noindent\textit{2000~Mathematics Subject Classification:} 35L65.

  \medskip


\end{abstract}

\Section{Introduction and Main Result}
\label{sec:Intro}

Consider the following nonlinear system of balance laws:
\begin{equation}
  \label{eq:BL}
  \partial_t u + \partial_x f(u) = G(u)
\end{equation}
where $f$ is the flow of a nonlinear hyperbolic system of conservation
laws and $G \colon \L1 \mapsto \L1$ is a (possibly) \emph{non local}
operator. It is known, see~\cite[Theorem~2.1]{ColomboGuerra}, that for
small times equation~\Ref{eq:BL} generates a Lipschitz semigroup.

When the source term is \emph{dissipative}, the existence of solutions
can be proved for all times, see~\cite{DafermosJHDE, DafermosHsiao} as
well as the continuous dependence, see~\cite{AmadoriGuerra2002,
  Cleo2}. These papers all deal with local sources. Memory effects,
i.e.~sources non local in time, were recently considered, for
instance, in~\cite{Cleo1}.

Here, we deal with dissipative non local sources and we provide the
well posedness of the Cauchy problem for~\Ref{eq:BL}, as well as
estimates on the dependence of the solutions from $f$ and $G$. The
proof relies on the combination of the Standard Riemann Semigroup $S$
generated by the conservation law $\partial_t u + \partial_x f(u) =
0$, see~\cite[Definition~9.1]{BressanLectureNotes}, combined through
the operator splitting technique with the Euler polygonal $(t,u)
\mapsto u + t \, G(u)$ generated by the ordinary differential system
$\partial_t u = G(u)$.

Here, we limit our attention to right hand sides of the type
\begin{equation}
  \label{eq:source}
  G(u) = g(u) + Q * u
\end{equation}
where $g \in \C{1,1}(\Omega;\reali^n)$ and $Q \in
\L1(\reali;\reali^n)$, the convolution being in the space variable,
see~\cite[\S~2]{ColomboGuerra} for several physical motivations.

Let $R$ be the matrix whose columns are the right eigenvectors of
$Df(0)$. We call the source term~\Ref{eq:source} \emph{column
  diagonally dominant}, see~\cite{DafermosHsiao}, if there exists a $c
> 0$ such that for $i=1, \ldots, n$ the matrix $M = R^{-1} Dg(0) \, R$
satisfies
\begin{equation}
  \label{eq:c}
  M_{ii} + \sum_{j=1, \, j \neq i}^n \modulo{M_{ji}} < -c\,,
\end{equation}
see~\cite[formula~(5)]{AmadoriGuerra1999} for a coordinate independent
extension of diagonal dominance.

It is well known that this dissipativity condition allows to prove the
well posedness globally in time of the Cauchy problem
for~\Ref{eq:BL}--\Ref{eq:source} in the case $Q = 0$ and $g \colon
\reali^n \mapsto \reali^n$, see~\cite{AmadoriGuerra2002,
  DafermosHsiao}. Similar global results can be obtained by means of
suitable $\L1$ estimates for relevant classes of systems,
see~\cite{DafermosJHDE}.

Our main result is the following.

\begin{theorem}
  \label{thm:BL}
  Fix an open set $\Omega \subseteq \reali^n$, with $0 \in \Omega$,
  and assume that
  \begin{description}
  \item[(F)] $f \in \C4(\Omega;\reali^n)$ is such that $Df$ is
    strictly hyperbolic with each characteristic field either
    genuinely nonlinear or linearly degenerate;
  \item[(G)] $g \in \C{1,1}(\Omega;\reali^n)$, $g(0)=0$, $g$ is column
    diagonally dominant and $Q \in \L1(\reali;\reali^n)$.
  \end{description}
  \noindent Then, there exist positive $\delta_1,\delta_2,
  \mathcal{L}, \kappa$ such that for all $Q$ with $\norma{Q}_{\L1}
  \leq \delta_1$ there exists a global semigroup $P \colon \left[0,
    +\infty\right[ \times \mathcal{D} \mapsto \mathcal{D}$ such that
  \begin{enumerate}[(1)]
  \item $\mathcal{D} \supseteq \left\{ u \in \L1 (\reali;\Omega)
      \colon \tv(u) \leq \delta_2 \right\}$;
  \item \label{it:3} for all $u_o \in \mathcal{D}$ and for $t \in
    \left[0, +\infty \right[$, the map $(t,x) \mapsto (P_{t} u_o)(x)$
    is a weak entropy solution to~\Ref{eq:BL} with initial datum
    $u_o$;
  \item for $t,s \in \left[0, +\infty \right[$, $u,w \in \mathcal{D}$
    and $s < t$, then
    \begin{equation}
      \label{eq:error}
      \begin{array}{rcl}
        \displaystyle
        \norma{P_t u - P_t w}_{\L1} 
        & \leq &
        \displaystyle 
        \mathcal{L}\cdot e^{-\kappa t}_{\phantom{\vert}} \cdot \norma{u-w}_{\L1}
        \\
        \displaystyle
        \norma{P_t u - P_s u}_{\L1} 
        & \leq &
        \displaystyle
        \mathcal{L} \cdot \left( 1 + \norma{u}_{\L1} \right)^{\phantom{\vert}}
        \cdot \modulo{t-s}\,;
      \end{array}
    \end{equation}
  \item \label{it:tg} if $S$ is the SRS generated by $\partial_t u +
    \partial_x f(u) = 0$, then for all $u \in \mathcal{D}$,
    \begin{displaymath}
      \lim_{t \to 0}
      \frac{1}{t} \,
      \norma{P_t u - \left( S_t u + t \, G(u) \right) }_{\L1} =0 \,;
    \end{displaymath}
  \end{enumerate}
  \noindent Moreover, let $f, \tilde f$ both satisfy~\textbf{(F)},
  both pairs $g,Q$ and $\tilde g, \tilde Q$
  satisfy~\textbf{(G)}. Denote by $P, \tilde P$ the corresponding
  processes and $\mathcal{D}_\delta, \tilde\mathcal{D}_{\tilde\delta}$
  their domains. Choose $\delta, \tilde \delta$ so that $\tilde
  \mathcal{D}_{\tilde \delta} \subseteq \mathcal{D}_\delta$. Then, for
  all $u \in \tilde \mathcal{D}_{\tilde \delta}$
  \begin{equation}
    \label{eq:Uffa}
    \begin{array}{rcl}
      \displaystyle
      \norma{P_t u - \tilde P_t u}_{\L1}
      & \leq &
      \displaystyle
      \mathcal{L} \cdot 
      \norma{Df - D \tilde f_2}_{\C0(\Omega;\reali^{n\times n})}
      \cdot t
      \\[4pt]
      & &
      \displaystyle
      + \, 
      \mathcal{L} \cdot
      \norma{g - \tilde g}_{\C0(\Omega;\reali^n)}
      \cdot t
      \\[4pt]
      & &
      \displaystyle
      + \, 
      \mathcal{L} \cdot
      \norma{Q - \tilde Q}_{\L1(\reali;\reali^n)}
      \cdot t \,.
    \end{array}
  \end{equation}
\end{theorem}

Condition~\Ref{it:tg} ensures that the orbits of $P$ are weak entropy
solutions, see~\cite[Corollary~3.13]{ColomboGuerra}. Moreover, the
solution yielded by $P$ can be characterized also as \emph{viscosity
  solution} in the sense of the integral inequalities
in~\cite[\S~9.2]{BressanLectureNotes}, see~\cite[(6) and~(7) in
Theorem~1.2]{ColomboGuerra}. Theorem~\ref{thm:BL} is obtained
applying~\cite[Theorem~2.5]{ColomboGuerra3} with $X =
\L1(\reali;\bar\Omega)$, see Section~\ref{sec:Tech}.

\Section{Outline of the Proof}
\label{sec:Sketch}

We sketch below the procedure used to prove Theorem~\ref{thm:BL}. All
technical details are deferred to Section~\ref{sec:Tech}.

Our general reference for the basic notions related to systems of
conservation laws is~\cite{BressanLectureNotes}. We assume throughout
that $0 \in \Omega$ and that $f$ satisfies~\textbf{(F)} in
Theorem~\ref{thm:BL}. Let $\lambda_1(u), \ldots, \lambda_n(u)$ be the
$n$ real distinct eigenvalues of $D\!f(u)$, indexed so that $\lambda_j
(u) < \lambda_{j+1}(u)$ for all $j$ and $u$. The $j$-th right,
respectively left, eigenvector is denoted $r_j(u)$, respectively
$l_j(u)$.

Let $\sigma \mapsto R_j(\sigma)(u)$, respectively $\sigma \mapsto
S_j(\sigma)(u)$, be the $j$-rarefaction curve, respectively the
$j$-shock curve, exiting $u$. If the $j$-th field is linearly
degenerate, then the parameter $\sigma$ above is the arc-length. In
the genuinely nonlinear case,
see~\cite[Definition~5.2]{BressanLectureNotes}, we choose $\sigma$ so
that
\begin{displaymath}
  \frac{\partial \lambda_j}{\partial\sigma}
  \left(R_j (\sigma)(u) \right) 
  =
  k_j
  \quad \mbox{ and } \quad
  \frac{\partial \lambda_j}{\partial\sigma}
  \left(S_j (\sigma)(u) \right) 
  =
  k_j \,,
\end{displaymath}
where $k_1,\ldots,k_n$ are positive and such that, as
in~\cite{AmadoriGuerra2002},
\begin{displaymath}
  \frac{\partial}{\partial\sigma} \left( R_j(\sigma)(0)\right) = r_j(0)
  \,, \quad
  \norma{r_j(0)} =1 \,.
\end{displaymath}
Introduce the $j$-Lax curve
\begin{displaymath}
  \sigma \mapsto \psi_j (\sigma) (u) =
  \left\{
    \begin{array}{c@{\qquad\mbox{ if }\quad}rcl}
      R_j(\sigma)(u) & \sigma & \geq & 0
      \\
      S_j(\sigma)(u) & \sigma & < & 0
    \end{array}
  \right.
\end{displaymath}
and for $\boldsymbol{\sigma} \equiv (\sigma_1, \ldots, \sigma_n)$,
define the map
\begin{displaymath}
  \mathbf{\Psi}(\boldsymbol{\sigma})(u^-)
  =
  \psi_n(\sigma_n)\circ\ldots\circ\psi_1(\sigma_1)(u^-) \,.
\end{displaymath}
By~\cite[\S~5.3]{BressanLectureNotes}, given any two states $u^-,u^+
\in \Omega$ sufficiently close to $0$, there exists a map $E$ such
that
\begin{equation}
  \label{eq:E}
  \boldsymbol{\sigma} = E(u^-,u^+)
  \quad \mbox{ if and only if } \quad
  u^+ = \mathbf{\Psi}(\boldsymbol{\sigma})(u^-) \,.
\end{equation}
elementary computations show that
\begin{equation}
  \label{eq:DE}
  D_u E(0,u)_{\Big\vert u=0} 
  = 
  R^{-1} 
  = 
  \left[
    \begin{array}{c}
      l_1(0) \\ \vdots \\ l_n(0)
    \end{array}
  \right] \,.
\end{equation}
Similarly, let the map $\mathbf{S}$ and the vector $\mathbf{q} = (
q_1, \ldots, q_n)$ be defined by
\begin{equation}
  \label{eq:S}
  u^+
  =
  \mathbf{S}(\mathbf{q})(u^-) =
  S_n(q_n) \circ \ldots \circ S_1(q_1) (u^-)
\end{equation}
as the gluing of the Rankine - Hugoniot curves.

Let $u$ be piecewise constant with finitely many jumps and assume that
$\tv(u)$ is sufficiently small. Call $\mathcal{I}(u)$ the finite set
of points where $u$ has a jump.  Let $\sigma_{x,i}$ be the strength of
the $i$-th wave in the solution of the Riemann problem for
\begin{equation}
  \label{eq:HCL}
  \partial_t u + \partial_x f(u) =0
\end{equation}
with data $u(x-)$ and $u(x+)$, i.e.~$(\sigma_{x,1}, \ldots,
\sigma_{x,n}) = E\left( u(x-), u(x +) \right)$. Obviously if $x\not\in
\mathcal{I}(u)$ then $\sigma_{x,i}=0$, for all $i=1,\ldots,n$.  As
in~\cite[\S~7.7]{BressanLectureNotes}, $\mathcal{A}(u)$ denotes the
set of approaching waves in $u$:
\begin{displaymath}
  \mathcal{A} (u) 
  = \left\{
    \begin{array}{c}
      \left(
        (x,i),(y,j)\right) \in \left( \mathcal{I}(u) \times \{1,\ldots,n\} 
      \right)^2
      \colon \\
      x < y \mbox{ and either } i > j \mbox{ or } i = j, 
      \mbox{ the $i$-th field}\\
      \mbox{is genuinely non linear, }
      \min \left\{ \sigma_{x,i}, \sigma_{y,j} \right\} <0\!  
    \end{array}
  \right\}
\end{displaymath}
while the linear and the interaction potential, see~\cite{Glimm}
or~\cite[formula~(7.99)]{BressanLectureNotes}, are
\begin{displaymath}
  \mathbf{V}(u)
  =
  \sum_{x\in I(u)} \sum_{i=1}^n
  \modulo{\sigma_{x,i}}
  \quad \mbox{ and } \quad
  \mathbf{Q}(u)
  =
  \sum_{\left((x,i),(y,j)\right) \in \mathcal{A}(u)}
  \modulo{\sigma_{x,i}\sigma_{y,j}} \,.
\end{displaymath}
Moreover, let
\begin{equation}
  \label{def:ups}
  \mathbf{\Upsilon} (u) 
  = 
  \mathbf{V} (u) + C_0 \cdot \mathbf{Q} (u)  
\end{equation}
where $C_0>0$ is the constant appearing in the functional of the
wave--front tracking algorithm,
see~\cite[Proposition~7.1]{BressanLectureNotes}. Recall that $C_0$
depends only on the flow $f$ and the upper bound of the total
variation of initial data.

Finally we define
\begin{eqnarray}
  \label{def:2.6}
  \mathcal{D}_\delta^*
  & = &
  \left\{
    v \in \L1\left(\reali,\Omega \right) \colon
    v \hbox { piecewise constant and } \mathbf{\Upsilon}(v) < \delta
  \right\}
  \\
  \nonumber
  \mathcal{D}_\delta
  & = &
  \mathrm{cl}\left\{\mathcal{D}_\delta^*\right\}
\end{eqnarray}
where the closure is in the strong $\L1$--topology. Observe that
$\mathcal{D}_\delta$ contains all $\L1$ functions with sufficiently
small total variation.

We now pass to the stability functional introduced
in~\cite{BressanYangLiu, LiuYang1, LiuYang3}. For any $\bar
v\in\mathcal{D}_\delta^*$, denote by $\bar\sigma_{x,i}$ the size of
the $i$--wave in the solution of the Riemann Problem with data $\bar
v(x-)$ and $\bar v(x+)$. Then define
\begin{displaymath}
  A_j^- [\bar v](x)
  =
  \sum_{y\leq x} \modulo{\bar\sigma_{y,j}},
  \quad 
  A_j^+ [\bar v](x)
  =
  \sum_{y> x} \modulo{\bar\sigma_{y,j}},
  \quad
  \mbox{ for } j = 1,\ldots,n\,.
\end{displaymath}
If the $i$-th characteristic field is linearly degenerate, then define
$\mathbf{A}_i$ as
\begin{equation}
  \label{eq:LinDeg}
  \mathbf{A}_i\left[\bar v\right]\left(q,x\right)
  =
  \sum_{1\leq j<i}A_j^+\left[\bar v\right](x) + 
  \sum_{ i<j\leq n}A_j^-\left[\bar v\right](x).
\end{equation}
While if the $i$-th characteristic field is genuinely nonlinear
\begin{equation}
  \label{eq:BoldA}
  \begin{array}{rcl}
    \mathbf{A}_i [\bar v] (q,x)
    & = &
    \displaystyle
    \sum_{1\leq j<i} A_j^+ [\bar v] (x) 
    + 
    \sum_{ i<j\leq n}A_j^- [\bar v] (x)
    \\[16pt]
    & &
    +
    A_i^+ [\bar v] (x) \cdot \caratt{\left[0, +\infty\right[} (q)
    +
    A_i^- [\bar v] (x) \cdot \caratt{\left]-\infty, 0\right[} (q).
  \end{array}
\end{equation}
Now choose $v,\tilde v$ piecewise constant in $\mathcal{D}_\delta^*$
and define the weights
\begin{equation}
  \label{def:w}
  \begin{array}{rcl}
    \mathbf{W}_i [v,\tilde v](q,x)
    & = &
    1
    +
    \kappa_1 \mathbf{A}_i [v] (q,x)
    +
    \kappa_1 \mathbf{A}_i [\tilde v] (-q,x)
    \\[2pt]
    & &
    \quad +
    \kappa_1 \kappa_2 \left(\mathbf{Q} (v) + \mathbf{Q} (\tilde v) \right)\,.
  \end{array}
\end{equation}
the constants $\kappa_1$ and $\kappa_2$ being those
in~\cite[Chapter~8]{BressanLectureNotes}. Define implicitly the
function $\mathbf{q}(x) \equiv \left( q_1(x), \ldots, q_n(x) \right)$
by
\begin{displaymath}
  \tilde v(x)
  =
  \mathbf{S} \left( \mathbf{q}(x) \right) \left(v(x)\right)
\end{displaymath}
with $\mathbf{S}$ as in~\Ref{eq:S}. The stability functional $\gf$ is
\begin{equation}
  \label{eq:Phi}
  \gf(v,\tilde v)
  =
  \sum_{i=1}^n \int_{-\infty}^{+\infty}
  \modulo{q_i(x)} \cdot 
  \mathbf{W}_i [v,\tilde v] \left(q_i(x),x\right) \, dx.
\end{equation}

We stress that $\gf$ is slightly different from the functional $\Phi$
defined in~\cite[formula~(8.6)]{BressanLectureNotes}. Indeed, here
\emph{all} jumps in $v$ or in $\tilde v$ are considered. There, on the
contrary, exploiting the structure of $\epsilon$-approximate front
tracking solutions, see~\cite[Definition~7.1]{BressanLectureNotes}, in
the definition of $\Phi$ the jumps due to non physical waves are
neglected when defining the weights $A_i$ and are considered as
belonging to a fictitious $(n+1)$-th family in the
definition~\cite[formula~(7.54)]{BressanLectureNotes} of $Q$.

Recall the following basic result in the theory of non linear systems
of conservation laws.

\begin{theorem}
  \label{thm:SRS}
  Let $f$ satisfy~\textbf{(F)}. Then, there exists a positive
  $\delta_o$ such that the equation~\Ref{eq:HCL} generates for all
  $\delta \in \left]0, \delta_o \right[$ a Standard Riemann Semigroup
  (SRS) $S\colon \left[0, +\infty \right[ \times \mathcal{D}_\delta
  \mapsto \mathcal{D}_\delta$, with Lipschitz constant $L$.
\end{theorem}

We refer to~\cite[Chapters~7 and~8]{BressanLectureNotes} for the proof
of the above result as well as for the definition and further
properties of the SRS.

Recall the following result from~\cite{ColomboGuerra2}:

\begin{proposition}
  \label{prop:Functionals}
  The functionals $\mathbf{\Upsilon}$, $\mathbf{Q}$ and
  $\mathbf{\Phi}$ admit an $\L1$ lower semicontinuous extension to all
  $\mathcal{D}_\delta$. Moreover,
  \begin{enumerate}
  \item for all $u \in \mathcal{D}_\delta$, the maps $t \mapsto
    \mathbf{Q} (S_t u)$ and $t \mapsto \mathbf{\Upsilon} (S_t u)$ are
    non increasing.
  \item for all $u,v \in \mathcal{D}_\delta$, the map $t \mapsto
    \mathbf{\Phi}(S_tu, S_t v)$ is non increasing;
  \item there exists a positive $C$ such that for all $u \in
    \mathcal{D}_\delta$,
    \begin{displaymath}
      \frac{1}{C} \tv(u) 
      \leq 
      \mathbf{\Upsilon} (u) 
      \leq
      C\, \tv(u)
    \end{displaymath}
    and for all $u,v \in \mathcal{D}_\delta$,
    \begin{displaymath}
      \frac{1}{C} \norma{u - \tilde u}_{\L1} 
      \leq 
      \mathbf{\Phi}(u,\tilde u) 
      \leq 
      C \, \norma{u - \tilde u}_{\L1} \,;
    \end{displaymath}
  \item for all $u \in \mathcal{D}_\delta$,
    \begin{eqnarray*}
      \mathbf{Q}(u) 
      & = &
      \liminf_{v \in \mathcal{D}^*_\delta, v \to u} \mathbf{Q}(v)
      \\
      \mathbf{\Upsilon}(u) 
      & = &
      \liminf_{v  \in \mathcal{D}^*_\delta, v \to u}
      \mathbf{\Upsilon}(v)
      \\
      \mathbf{\Phi}(u, \tilde u)
      & = &
      \liminf_{v \in \mathcal{D}^*_\delta, v \to u} \mathbf{\Phi}(v, \tilde v)
    \end{eqnarray*}
  \end{enumerate}
\end{proposition}

\noindent The results in~\cite{ColomboGuerra2} also provide an
explicit expression of $\mathbf{\Phi}$ in terms of \emph{wave
  measures}, see~\cite[\S~10.1]{BressanLectureNotes}. For the
properties of $\mathbf{Q}$ and $\mathbf{\Upsilon}$, see
also~\cite{BaitiBressan2, BressanLectureNotes,BressanColombo2}.

Introduce the map
\begin{eqnarray}
  \label{eq:F}
  \F(s)u
  =
  u + s \, g(u) + s \, Q*u .
\end{eqnarray}
that satisfies the properties stated in the following lemma, whose
proof is deferred to Section~\ref{sec:Tech}.

\begin{lemma}
  \label{lem:Stime}
  Let~\textbf{(F)} and~\textbf{(G)} hold. For all $\delta$
  sufficiently small and all $u,\tilde u \in \mathcal{D}_\delta$,
  \begin{eqnarray*}
    \mathbf{Q} \left( \F(s)u \right)
    & \leq &
    \mathbf{Q}(u) + \O s \, \mathbf{\Upsilon}^2 (u)
    \\
    \mathbf{\Upsilon} \left( \F(s)u \right)
    & \leq &
    \left( 1 - \frac{c}{8} s \right) \mathbf{\Upsilon}(u)
    \\
    \mathbf{\Phi} \left( \F(s)u , \F(s)\tilde u \right)
    & \leq &
    \left( 1- \frac{c}{4}s \right) \mathbf{\Phi}(u, \tilde u) \,.
  \end{eqnarray*}
\end{lemma}

We now recall the basic definitions and results
from~\cite[Section~2]{ColomboGuerra3} that allow us to complete the
proof of Theorem~\ref{thm:BL}. In the complete metric space $X =
\L1(\reali;\bar\Omega)$ with the $\L1$ distance, select for a fixed
$M$ the closed domain
\begin{displaymath}
  \mathcal{D}^M
  = 
  \left\{
    u \in \mathcal{D}_\delta \colon \mathbf{\Phi}(u,0) \leq M
  \right\},
\end{displaymath}
and the local flow
\begin{equation}
  \label{eq:LocalFlow}
  F \colon [0,\tau]\times \mathcal{D}^M \mapsto \mathcal{D}^M
  \qquad \qquad
  F(s)u = \F(s) S_s u \,,
\end{equation}
where $S$ is the SRS of Theorem~\ref{thm:SRS}, $\tau$ is positive and
sufficiently small. Recall that, in the present autonomous setting, a
Lipschitz continuous map is a \emph{local flow}
by~\cite[Definition~2.1]{ColomboGuerra3}.  Note that $\mathcal{D}^M$
is invariant with respect to both $\F$ and $S$, so that the above
definition makes sense.

For any positive $\epsilon$, the Euler $\epsilon$-polygonal generated
by $F$ is
\begin{equation}
  \label{eq:polygonal}
  F^\epsilon (t) \, u
  =
  F(t-k\epsilon) \circ \comp_{h=0}^{k-1} F(\epsilon) \, u
\end{equation}
see~\cite[Definition~2.2]{ColomboGuerra3}. Therein, the following
theorem is proved in a generic complete metric space.

\begin{theorem}
  \label{thm:main}
  Let $F\colon [0,\tau]\times \hat \mathcal{D} \mapsto \hat
  \mathcal{D}$ be a local flow that satisfies
  \begin{enumerate}
  \item \label{it:first} there exists a non decreasing map $\omega
    \colon [0,\tau/2] \mapsto \reali^+$ with $\int_0^{\tau/2}
    \frac{\omega(\xi)}{\xi} \, d\xi < +\infty$ such that for all
    $(s,u)$ and all $k \in \naturali$
    \begin{equation}
      \label{eq:k}
      d \left( F (k s) \circ F(s) u, F \left( (k+1) s \right) u
      \right)
      \leq
      k s \, \omega(s) \,;
    \end{equation}
  \item \label{it:second} there exists a positive $L$ such that for
    all $\epsilon \in [0, \tau]$ and for all $t \geq 0$
    \begin{equation}
      \label{eq:Stability}
      d \left( F^\epsilon (t) u, F^\epsilon (t) \tilde u \right)
      \leq L \cdot d(u,\tilde u) \,.
    \end{equation}
  \end{enumerate}
  \noindent Then, there exists a unique Lipschitz semigroup $P \colon
  \left[0, +\infty\right[ \times \hat \mathcal{D} \mapsto \hat
  \mathcal{D}$ such that for all $u \in \hat \mathcal{D}$
  \begin{equation}
    \label{eq:tangent}
    \frac{1}{s} \,
    d\left( P_su, F(s)u \right)
    \leq
    \frac{2L}{\ln 2} \cdot 
    \int_0^s \frac{\omega(\xi)}{\xi} \, d\xi \,.
  \end{equation}
\end{theorem}
The proof of Theorem~\ref{thm:BL} is deferred to
Section~\ref{sec:Tech}. It amounts to show that the above abstract
result can be applied in the present setting, with $\hat \mathcal{D} =
\mathcal{D}^M$ and $F$ as in~\Ref{eq:LocalFlow}.

\Section{Technical Details}
\label{sec:Tech}

\begin{lemma}
  \label{lem:Taylor}
  Let $f$ satisfy~\textbf{(F)}, $\Omega$ be a sufficiently small
  neighborhood of the origin; $a,b \in \reali^n$ and $s \geq 0$ be
  sufficiently small. Choose $u^-,v^- \in \Omega$ and define $u^+ =
  u^- + s a $ and $v^+ = v^- + s b$. Then, if $\boldsymbol{\sigma^-}$,
  $\boldsymbol{\sigma^+}$ satisfy $v^- = \mathbf{\Psi}
  (\boldsymbol{\sigma^-}) (u^-)$ and $v^+ = \mathbf{\Psi}
  (\boldsymbol{\sigma^+}) (u^+)$,
  \begin{equation}
    \label{SizeEstimate}
    \sum_{i=1}^n \modulo{\sigma^+_i - \sigma^-_i}
    \leq
    \O \cdot
    s
    \left( \sum_{i=1}^n \modulo{\sigma_i^-} + \norma{b-a} \right) \,.
  \end{equation}
  If $g$ satisfies~\textbf{(G)}, $u^+ = u^- + s \left(a + g(u^-)
  \right)$ and $v^+ = v^- + s \left(b + g(v^-) \right)$, then
  \begin{equation}
    \label{SizeEstimateBIS}
    \sum_{i=1}^n \modulo{\sigma^+_i}
    \leq
    \left(1 - \frac{c}{2} s \right) \sum_{i=1}^n \modulo{\sigma_i^-} 
    + \O \, s \, \norma{b-a} \,.
  \end{equation}
  An entirely analogous result holds with the map $\mathbf{\Psi}$
  replaced by the gluing $\mathbf{S}$ of shock curves, i.e.~$v^- =
  \mathbf{S} (\boldsymbol{\sigma^-}) (u^-)$ and $v^+ = \mathbf{S}
  (\boldsymbol{\sigma^+}) (u^+)$.
\end{lemma}

\begin{proof}
  Let $\boldsymbol{\sigma}^- = \boldsymbol{\sigma}$ and
  $\boldsymbol{\sigma}^+ = E \left( u + s a, \mathbf{\Psi}
    (\boldsymbol{\sigma})(u) + s b \right)$. Introduce the $\C2$ map
  \begin{displaymath}
    \phi(s,a,b,\boldsymbol{\sigma}) 
    = 
    E \left(u+ s a,\mathbf{\Psi}(\boldsymbol{\sigma})(u) + s b\right) 
    - 
    \boldsymbol{\sigma} \,.
  \end{displaymath}
  Note that $\phi(0,a,b,\sigma) = 0$ and $\phi(s,a,a,0) = 0$,
  by~\cite[Lemma~2.5]{BressanLectureNotes} we get
  \begin{displaymath}
    \norma{\phi(s,a,b,\boldsymbol{\sigma})}
    \leq 
    \O
    s 
    \left( \sum_{i=1}^n \modulo{\sigma_i} + \norma{b-a} \right)
  \end{displaymath}
  proving~\Ref{SizeEstimate}.

  To prove~\Ref{SizeEstimateBIS}, introduce the functions
  \begin{displaymath}
    B_{ij}(a,u)
    = 
    \frac{\partial^2}{\partial s \, \partial \sigma_j}
    E_i \left( 
      \begin{array}{c}
        u + s \left( a+ g(u) \right), 
        \\
        \mathbf{\Psi}(\boldsymbol{\sigma})(u) + s 
        \left( 
          b + g \left( \mathbf{\Psi}(\boldsymbol{\sigma})(u)\right) 
        \right)
      \end{array}
    \right)_{\strut \left\vert
        \begin{array}{l}
          \scriptstyle \sigma=0, \\[-4pt]
          \scriptstyle s=0,\\[-4pt] 
          \scriptstyle b=a, 
        \end{array}
      \right.
    }
  \end{displaymath}
  \begin{displaymath}
    \phi(s,a,b,\sigma,u)
    =
    E \left( \!
      \begin{array}{c}
        u + s \left( a+ g(u) \right), 
        \\
        \mathbf{\Psi}(\boldsymbol{\sigma})(u) + s 
        \left( 
          b + g \left( \mathbf{\Psi}(\boldsymbol{\sigma})(u)\right) 
        \right)
      \end{array}
      \! \right)
    - \left( \Id + s B(a,u) \right) \boldsymbol{\sigma} .
  \end{displaymath}
  By the $\C{2,1}$ regularity of $E$, the $n \times n$ matrix $B$ is a
  Lipschitz function of $(a,u)$. Moreover, by~\Ref{eq:DE}
  \begin{eqnarray*}
    B_{ij}(0,0)
    & = &
    \frac{\partial^2}{\partial s \, \partial \sigma_j}
    E_i \left( 
      0, 
      \psi_j(\sigma_j)(0) + s g \left( \psi_j(\sigma_j)(0) \right)
    \right)_{\strut \Big\vert s=0, \sigma_j=0}
    \\
    & = &
    \frac{\partial}{\partial s}
    \left(
      l_i(0) \cdot \left( r_i(0) + s \, D g(0) \, r_j(0) \right)
    \right)_{\Big\vert s=0}
    \\
    & = &
    l_i(0) \, D g(0) \, r_j(0)
    \\
    & = &
    M_{ij} \,.
  \end{eqnarray*}
  Therefore, by~\Ref{eq:c}
  \begin{equation}
    \label{eq:Estimate}
    \sum_{i=1}^n \modulo{ \left[ 
        \left( \Id + s B(0,0) \right) \boldsymbol{\sigma}
      \right]_i}
    \leq
    (1-cs) \sum_{i=1}^n \modulo{\sigma_i} \,.
  \end{equation}

  Note that $\phi(0,a,b,\boldsymbol{\sigma},u)=0$, hence by the
  Lipschitzeanity of $D\phi$,
  \begin{displaymath}
    \norma{\phi(s,a,b,\boldsymbol{\sigma},u) - 
      \phi(s,a,a,\boldsymbol{\sigma},u)}
    \leq
    \O s \norma{b-a} \,.
  \end{displaymath}
  We thus have
  \begin{eqnarray}
    \nonumber
    \norma{\phi(s,a,b,\boldsymbol{\sigma},u)}
    & \leq &
    \norma{\phi(s,a,b,\boldsymbol{\sigma},u) - 
      \phi(s,a,a,\boldsymbol{\sigma},u)}
    +
    \norma{\phi(s,a,a,\boldsymbol{\sigma},u)}
    \\
    \label{eq:phi}
    & \leq &
    \O s \norma{b-a}
    +
    \norma{\phi(s,a,a,\boldsymbol{\sigma},u)} \,.
  \end{eqnarray}
  By the definition of $\phi$ and the choice of $B_{ij}(a,u)$,
  \begin{eqnarray*}
    \phi(0,a,a,\boldsymbol{\sigma}, u)
    & = &
    0
    \\
    \phi(s,a,a,0,u)
    & = &
    0
    \\
    \frac{\partial^2}{\partial s \, \partial \sigma_j} \,
    \phi(s,a,b,\boldsymbol{\sigma},u)_{\left\vert
        \begin{array}{l}
          \scriptstyle \boldsymbol{\sigma}=0, \\ [-4pt]
          \scriptstyle s=0, \\[-4pt]
          \scriptstyle b=a
        \end{array}\right.
    }
    & = &
    0
  \end{eqnarray*}
  using again~\cite[Lemma~2.5]{BressanLectureNotes}, we obtain
  \begin{equation}
    \label{eq:phiphi}
    \norma{\phi(s,a,a,\boldsymbol{\sigma},u)}
    \leq
    \O s \sum_{i=1}^n \modulo{\sigma_i} 
    \left( s + \sum_{i=1}^n \modulo{\sigma_i}\right) \,.
  \end{equation}
  Finally,
  \begin{eqnarray*}
    & &
    E \left( 
      u + s \left( a+ g(u) \right), 
      \mathbf{\Psi}(\boldsymbol{\sigma})(u) + s 
      \left( 
        b + g \left( \mathbf{\Psi}(\boldsymbol{\sigma})(u)\right) 
      \right)
    \right)
    \\
    & = &
    \phi(s,a,b,\boldsymbol{\sigma},u) 
    + 
    \left( \Id + s B(a,u) \right) \boldsymbol{\sigma}
    \\
    & = &
    \left( \Id + s B(0,0) \right) \boldsymbol{\sigma}
    +
    \phi(s,a,b,\boldsymbol{\sigma},u) 
    + 
    s \left( B(a,u) - B(0,0) \right) \boldsymbol{\sigma} \,.
  \end{eqnarray*}
  Apply now~\Ref{eq:Estimate}, \Ref{eq:phi}, \Ref{eq:phiphi} and the
  Lipschitzeanity of $B$
  \begin{eqnarray*}
    \sum_{i=1}^n \modulo{\sigma_i^+}
    & = &
    \sum_{i=1}^n 
    \modulo{  E_i \left( 
        \begin{array}{c}
          u + s \left( a+ g(u) \right), 
          \\
          \mathbf{\Psi}(\boldsymbol{\sigma})(u) + s 
          \left( 
            b + g \left( \mathbf{\Psi}(\boldsymbol{\sigma})(u)\right) 
          \right)
        \end{array}
      \right)}
    \\
    & \leq &
    \sum_{i=1}^n 
    \modulo{
      \left[ \left( \Id + s B(0,0) \right) \boldsymbol{\sigma} \right]_i
    }
    +
    \O s \, \norma{b-a}
    \\
    & &
    +
    \norma{\phi(s,a,a,\boldsymbol{\sigma},u)}
    +
    s\, \norma{B(a,u) - B(0,0)} \, \norma{\boldsymbol{\sigma}}
    \\
    & \leq &
    (1-cs) \sum_{i=1}^n \modulo{\sigma_i}
    + 
    \O s \norma{b-a}
    \\
    & &
    +
    \O s \sum_{i=1}^n \modulo{\sigma_i}  
    \left(s + \sum_{i=1}^n \modulo{\sigma_i} + \norma{a} +\norma{u} \right)
    \\
    & \leq &
    \left(1 - \frac{c}{2}s \right) \sum_{i=1}^n \modulo{\sigma_i}
    + 
    \O s \norma{b-a}
  \end{eqnarray*}
  provided $s$, $\sum_{i=1}^n \modulo{\sigma_i}$, $\norma{a}$ and
  $\norma{u}$ are sufficiently small.
\end{proof}

Introduce for $N \in \naturali$ the projection
\begin{displaymath}
  \Pi_N (u)
  =
  N \sum_{k=-1-N^2}^{-1+N^2} \int_{k/N}^{(k+1)/N} u(\xi) \, d\xi \;
  \chi_{\strut \left]k/N, (k+1)/N \right]}
\end{displaymath}
as in~\cite[\S~3.2]{ColomboGuerra}. Note that $\Pi_N u$ is piecewise
constant. For later use, we introduce the following approximated local
flows, see~\cite[Definition~2.1]{ColomboGuerra3}, generated by the
source term
\begin{displaymath}
  \F_N(s)u
  =
  u + s \, g(u)+ s \, \Pi_N(Q*v).
\end{displaymath}

\begin{lemma}
  \label{lem:V}
  Let~\textbf{(F)} and~\textbf{(G)} hold. If $\delta$ and
  $\norma{Q}_{\L1}$ are sufficiently small, then for all $v \in
  \mathcal{D}_\delta^*$
  \begin{displaymath}
    \mathbf{V} \left( \F_N(s)v \right)
    \leq
    \left( 1 - \frac{c}{4} s \right) \mathbf{V}(v) \,.
  \end{displaymath}
\end{lemma}

\begin{proof}
  Denote $w = \Pi_N(Q*v)$ and introduce the piecewise constant
  function $v' = v + s g(v) + s w$. Let $\sigma_{x,i}$, respectively
  $\sigma'_{x,i}$ be the jumps in $v$, respectively $v'$.  Denote also
  $\Delta v(x) = v(x+)- v(x-)$ and similarly $\Delta w$, $\Delta v'$.

  Apply~\Ref{SizeEstimateBIS} at any $x \in \reali$, with $a = w(x-)$
  and $b = w(x+)$ to obtain
  \begin{eqnarray*}
    \mathbf{V}(v')
    & = &
    \sum_{x \in \reali} \sum_{i=1}^n \modulo{\sigma'_{x,i}}
    \\
    &\leq &
    \left( 1-\frac{c}{2} s \right) 
    \sum_{x\in\reali} \sum_{i=1}^n \modulo{\sigma_{x,i}}
    + \O s \, \sum_{x \in \reali} \norma{\Delta w(x)}
    \\
    & \leq &
    \left( 1-\frac{c}{2} s \right) \mathbf{V}(v)
    +
    \O s \tv \left( \Pi_N (Q*v)\right)
    \\
    & \leq &
    \left( 1-\frac{c}{2} s \right) \mathbf{V}(v)
    +
    \O s \tv (Q*v)
    \\
    & \leq &
    \left( 1-\frac{c}{2} s \right) \mathbf{V}(v)
    +
    \O s \norma{Q}_{\L1} \mathbf{V}(v)
    \\
    & \leq &
    \left( 1-\frac{c}{4} s \right) \mathbf{V}(v)
  \end{eqnarray*}
  provided $\norma{Q}_{\L1}$ is sufficiently small.
\end{proof}

\begin{proofof}{Lemma~\ref{lem:Stime}}
  By~\cite[Lemma~4.2]{ColomboGuerra2}, it is possible to choose a
  sequence $v^\nu \in \mathcal{D}^*_\delta$ such that $v^\nu \to u$ in
  $\L1$, $\mathbf{Q}(v^\nu)\to \mathbf{Q}(u)$ and
  $\mathbf{\Upsilon}(v^\nu) \to \mathbf{\Upsilon}(u)$.
  By~\cite[Proposition~1.1]{ColomboGuerra}, we may
  apply~\cite[Corollary~3.5]{ColomboGuerra} and use the
  estimate~\cite[(3.5) in Lemma~3.6]{ColomboGuerra} in the case
  $L_3=0$, $G(u) = g(u) + \Pi_N( Q*u)$. Note that the latter map is
  piecewise constant whenever $u$ is. We thus obtain
  \begin{eqnarray*}
    \mathbf{Q} \left( \F(s)u \right)
    & \leq &
    \liminf_{N \to +\infty}
    \mathbf{Q} \left( \F_N(s)u \right)
    \\
    & \leq &
    \liminf_{N \to +\infty}
    \liminf_{\nu \to +\infty}
    \mathbf{Q} \left( \F_N(s) v^\nu \right)
    \\
    & \leq &
    \liminf_{\nu \to +\infty}
    \left(
      \mathbf{Q}(v^\nu) + \O s \, \mathbf{V}^2(v^\nu)
    \right)
    \\
    & \leq &
    \liminf_{\nu \to +\infty}
    \left(
      \mathbf{Q}(v^\nu) + \O s \, \mathbf{\Upsilon}^2(v^\nu)
    \right)
    \\
    & = &
    \mathbf{Q}(u) + \O s \, \mathbf{\Upsilon}^2 (u)
  \end{eqnarray*}
  proving the former estimate.

  To prove the latter one, use the same sequence $v^\nu$,
  Lemma~\ref{lem:V} and follow an analogous argument based on the
  results in~\cite{ColomboGuerra}, to obtain
  \begin{eqnarray*}
    & &
    \mathbf{\Upsilon} \left( \F(s)u \right)
    \\
    & \leq &
    \liminf_{N \to +\infty} \liminf_{\nu \to +\infty}
    \mathbf{\Upsilon} 
    \left( \F_N(s) v^\nu \right)
    \\
    & \leq &
    \liminf_{\nu \to +\infty}
    \left(
      \left( 1 - \frac{c}{4}s \right) \mathbf{V}(v^\nu) 
      + 
      C_0 \mathbf{Q}(v^\nu)
      +
      \O s \mathbf{V}^2(v^\nu) 
    \right)
    \\
    & \leq &
    \liminf_{\nu \to +\infty}
    \left(
      \left( 1-\frac{c}{8}s \right) \mathbf{\Upsilon} (v^\nu)
      +
      \frac{c}{8}s 
      \left( 
        - \mathbf{V}(v^\nu) 
        + C_0 \mathbf{Q}(v^\nu)
        + \O \mathbf{V}^2(v^\nu)
      \right)
    \right)
    \\
    & \leq &
    \liminf_{\nu \to +\infty}
    \left(
      \left( 1-\frac{c}{8}s \right) \mathbf{\Upsilon} (v^\nu)
      +
      \frac{c}{8}s\, \mathbf{V}(v^\nu) 
      \left( 
        - 1 + \O \mathbf{V}(v^\nu)
      \right)
    \right)
    \\
    & \leq &
    \liminf_{\nu \to +\infty}
    \left( 1-\frac{c}{8}s \right) \mathbf{\Upsilon} (v^\nu)
    \\
    & \leq &
    \left( 1-\frac{c}{8}s \right) \mathbf{\Upsilon} (u)
  \end{eqnarray*}
  for a sufficiently small $\delta$.

  We now pass to the estimate on $\mathbf{\Phi}$.
  By~\cite[Lemma~4.5]{ColomboGuerra2}, we may choose two sequences of
  piecewise constant functions $v^\nu, \tilde v^\nu \in
  \mathcal{D}_\delta^*$ such that
  \begin{displaymath}
    v^\nu \to u 
    \quad \mbox{ and } \quad
    \tilde v^\nu \to \tilde u 
    \quad\mbox{ in }
    \L1,\qquad
    \mathbf{\Phi}(v^\nu, \tilde v^\nu) \to \mathbf{\Phi} (u, \tilde u) \,.
  \end{displaymath}
  Define implicitly the functions $\mathbf{q}^\pm(x)$ by $\F_N(s)
  (\tilde v^\nu) = \mathbf{S} (\mathbf{q}^+ ) \left( \F_N(s) v^\nu
  \right)$ and $\tilde v^\nu = \mathbf{S} (\mathbf{q}^-) (v^\nu)$.
  Moreover, with reference to~\Ref{def:w}, let
  \begin{eqnarray*}
    \mathbf{W}_i^+ (x)
    & = &
    \mathbf{W}_i 
    [\F_N(s) v^\nu, \F_N(s) \tilde v^\nu] 
    \left( \mathbf{q}_i^+(x), x \right)
    \\
    \mathbf{W}_i^- (x)
    & = &
    \mathbf{W}_i 
    [v^\nu, \tilde v^\nu] 
    \left( \mathbf{q}_i^-(x), x \right)
  \end{eqnarray*}
  and compute
  \begin{eqnarray}
    \nonumber
    & &
    \mathbf{\Phi} \left( \F_N(s) v^\nu, \F_N(s) \tilde v^\nu \right)
    -
    \mathbf{\Phi} (v^\nu, \tilde v^\nu)
    \\
    \nonumber
    & = &
    \int_{\reali} \sum_{i=1}^n
    \left(
      \modulo{q_i^+(x)} \mathbf{W}_i^+(x)
      -
      \modulo{q_i^-(x)} \mathbf{W}_i^-(x)
    \right)
    \, dx
    \\
    & = &
    \label{eq:R}
    \int_{\reali} 
    \mathcal{R}(x) \, dx
  \end{eqnarray}
  where the integrand $\mathcal{R}(x)$ is estimated splitting it as
  follows:
  \begin{eqnarray*}
    \mathcal{R}(x)
    & = &
    I(x) + I\!I(x) + \widetilde{I\!I}(x) + I\!I\!I(x) +
    \widetilde{I\!I\!I}(x)
    \\
    I(x)
    & = &
    \sum_{i=1}^n \modulo{q_i^+(x)} - \sum_{i=1}^n \modulo{q_i^-(x)}
    \\
    I\!I(x)
    & = &
    \kappa_1 
    \sum_{i=1}^n
    \Bigg( 
    \modulo{q_i^+(x)}
    \mathbf{A}_i[\F_N(s) v^\nu] \left( q_i^+(x),x \right)
    \\
    & &
    \qquad\qquad
    -
    \modulo{q_i^-(x)}
    \mathbf{A}_i[v^\nu] \left( q_i^-(x),x \right)
    \Bigg)
    \\
    \widetilde{I\!I}(x)
    & = &
    \kappa_1 
    \sum_{i=1}^n
    \Bigg( 
    \modulo{q_i^+(x)}
    \mathbf{A}_i[\F_N(s) \tilde v^\nu] \left( -q_i^+(x),x \right)
    \\
    & &
    \qquad\qquad
    -
    \modulo{q_i^-(x)}
    \mathbf{A}_i[\tilde v^\nu] \left( -q_i^-(x),x \right)
    \Bigg)
    \\
    I\!I\!I(x)
    & = &
    \kappa_1 \, \kappa_2 
    \left(
      \sum_{i=1}^n \modulo{q_i^+(x)} \mathbf{Q} \left( \F_N(s) v^\nu
      \right)
      -
      \sum_{i=1}^n \modulo{q_i^-(x)} \mathbf{Q} \left( v^\nu \right)
    \right)
    \\
    \widetilde{I\!I\!I}(x)
    & = &
    \kappa_1 \, \kappa_2 
    \left(
      \sum_{i=1}^n \modulo{q_i^+(x)} \mathbf{Q} \left( \F_N(s) \tilde v^\nu
      \right)
      -
      \sum_{i=1}^n \modulo{q_i^-(x)} \mathbf{Q} \left( \tilde v^\nu
      \right)
    \right)
  \end{eqnarray*}
  We now show that $\int_{\reali} I(x) \, dx$ is strictly negative and
  controls the growth in the other terms. By Lemma~\ref{lem:Taylor},
  applied to shock curves instead of Lax curves and with $a = a(x) =
  \left( \Pi_N Q*v^\nu \right) (x)$, $b = b(x) = \left(\Pi_N Q* \tilde
    v^\nu \right) (x)$
  \begin{displaymath}
    I(x)
    \leq
    -\frac{c}{2} s \sum_{i=1}^n \modulo{q_i(x)}
    +
    \O s \,
    \norma{b (x) - a (x)}
  \end{displaymath}
  Passing to the second addend, consider each term in the sum
  defining it:
  \begin{eqnarray}
    \nonumber
    & &
    \modulo{q_i^+(x)}
    \mathbf{A}_i[\F_N(s) v^\nu] \left( q_i^+(x),x \right)
    -
    \modulo{q_i^-(x)}
    \mathbf{A}_i[v^\nu] \left( q_i^-(x),x \right)
    \\
    \label{eq:Questa}
    & \leq &
    \modulo{q_i^+(x) - q_i^-(x)}
    \mathbf{A}_i[\F_N(s) v^\nu] \left( q_i^+(x),x \right)
    \\
    \label{eq:Quella}
    & &
    \quad
    +
    \modulo{q_i^-(x)}
    \modulo{
      \mathbf{A}_i[\F_N(s) v^\nu] \left( q_i^+(x),x \right)
      -
      \mathbf{A}_i[v^\nu] \left( q_i^-(x),x \right)
    }
  \end{eqnarray}
  By Lemma~\ref{lem:Taylor}, applied to shock curves,
  \begin{eqnarray}
    \nonumber
    & &
    \modulo{q_i^+(x) - q_i^-(x)}
    \\
    \nonumber
    & \leq &
    \O s 
    \Bigg( 
    \norma{g\left(v^\nu(x)\right) + \Pi_N Q * v^\nu (x)- 
      g\left(\tilde v^\nu(x)\right) - \Pi_N Q * \tilde v^\nu (x)} 
    \\
    \nonumber
    & &
    \qquad\qquad
    + 
    \sum_{i=1}^n \modulo{q_i^-(x)} \Bigg)
    \\
    \label{eq:qPiccola}
    & \leq &
    \O s \left( \sum_{i=1}^n \modulo{q_i^-(x)} + \norma{b(x) - a(x)} \right).
  \end{eqnarray}
  Note that $\mathbf{A}_i[\F_i(s) v^\nu](x) \leq \O\delta$, hence the
  first term in the right hand side of~\Ref{eq:Questa} is bounded by
  {\small
    \begin{displaymath}
      \modulo{q_i^+(x) - q_i^-(x)}
      \mathbf{A}_i[\F_N(s) v^\nu] \left( q_i^+(x),x \right)
      \leq
      \O s \delta \!
      \left( \sum_{i=1}^n \modulo{q_i^-(x)} + \norma{b(x) - a(x)} \right).
    \end{displaymath}
  }%
  Concerning the second term~\Ref{eq:Quella}, if $q_i^+(x) \cdot
  q_i^-(x) < 0$, then
  \begin{eqnarray*}
    \modulo{q_i^-(x)}
    & \leq &
    \modulo{q_i^+(x) - q_i^-(x)}
    \\
    & \leq &
    \O s \left( 
      \sum_{i=1}^n \modulo{q_i^-(x)} 
      + 
      \norma{b(x) - a(x)}
    \right)\,.
  \end{eqnarray*}
  Moreover, since $\modulo{\mathbf{A}_i[\F_i(s) v^\nu](x) -
    \mathbf{A}_i[v^\nu](x)} \leq \O \delta$, we get
  \begin{eqnarray*}
    & &
    \modulo{q_i^-(x)}
    \modulo{
      \mathbf{A}_i[\F_N(s) v^\nu] \left( q_i^+(x),x \right)
      -
      \mathbf{A}_i[v^\nu] \left( q_i^-(x),x \right)
    }
    \\
    & \leq &
    \O s\, \delta
    \left( 
      \sum_{i=1}^n \modulo{q_i^-(x)}  
      + 
      \norma{b(x) - a(x)}
    \right) \,.
  \end{eqnarray*}
  On the other hand, if $q_i^+(x) \cdot q_i^-(x) > 0$,
  by~\Ref{eq:LinDeg} or~\Ref{eq:BoldA}, a wave $\sigma'_{y,j}$ is
  counted in the sum defining $\mathbf{A}_i[\F_N(s)v^\nu]\left(
    q_i^+(x),x \right)$ if and only the wave $\sigma_{y,j}$ is counted
  in the sum defining $\mathbf{A}_i[v^\nu] \left( q_i^-(x),x
  \right)$. Therefore, by Lemma~\ref{lem:Taylor},
  \begin{eqnarray*}
    & &
    \modulo{
      \mathbf{A}_i[\F_N(s) v^\nu] \left( q_i^+(x),x \right)
      -
      \mathbf{A}_i[v^\nu] \left( -q_i^+(x),x \right)
    }
    \\
    & \leq &
    \sum_{x \in \reali} \sum_{i=1}^n 
    \modulo{\sigma'_{x,i} - \sigma_{x,i}}
    \\
    & \leq &
    \O s \sum_{x \in \reali}
    \left( 
      \sum_{i=1}^n 
      \modulo{\sigma_{x,i}}
      +
      \norma{\Delta a(x)}
    \right)
    \\
    & \leq &
    \O s \left( \mathbf{V}(v^\nu) + \tv ( \Pi_N Q * \tilde v^\nu )
    \right)
    \\
    & \leq &
    \O \, s \, \delta\,.
  \end{eqnarray*}

  The second addend $I\!I(x)$ is thus bounded as
  \begin{displaymath}
    I\!I(x)
    \leq
    \O \, s  \, \delta
    \left( 
      \sum_{i=1}^n \modulo{q_i^-(x)} + \norma{b(x) - a(x)}
    \right) \, .
  \end{displaymath}
  The term $\widetilde{I\!I}(x)$ can be treated repeating the same
  procedure

  Consider now each term in the sum defining $I\!I\!I(x)$, proceed as
  in~\Ref{eq:Questa}--\Ref{eq:Quella}, using~\Ref{eq:qPiccola} and
  similarly to the second part of Lemma~\ref{lem:Stime},
  \begin{eqnarray*}
    & &
    \left(
      \modulo{q_i^+(x)} \mathbf{Q} \left( \F_N(s) v^\nu
      \right)
      -
      \modulo{q_i^-(x)} \mathbf{Q} \left( v^\nu \right)
    \right)
    \\
    & \leq &
    \modulo{q_i^+(x) - q_i^-(x)} \mathbf{Q} \left(\F_N(s)v^\nu\right)
    +
    \modulo{q_i^-(x)} 
    \left( 
      \mathbf{Q}\left(\F_N(s)v^\nu\right) - \mathbf{Q}(v^\nu) 
    \right)
    \\
    & \leq &
    \O s \delta
    \left( \sum_{i=1}^n \modulo{q_i^-(x)} + \norma{b(x) - a(x)}\right)
    +
    \O s \delta \modulo{q_i^-(x)}
    \\
    & \leq &
    \O s \delta
    \left( \sum_{i=1}^n \modulo{q_i^-(x)} + \norma{b(x) - a(x)}\right).
  \end{eqnarray*}
  Obviously, $\widetilde{I\!I\!I}(x)$ is treated analogously. Thus,
  $\mathcal{R}$ in~\Ref{eq:R} is bounded as:
  \begin{displaymath}
    \mathcal{R}(x)
    \leq
    - \frac{c}{2} s \sum_{i=1}^n q_i(x) 
    + \O s \, \delta \sum_{i=1}^n \modulo{q_i(x)} 
    + \O s \, \norma{b(x) - a(x)} \,.
  \end{displaymath}
  Integrating and using the Lipschitzeanity of $\Pi_N$, a standard
  inequality on the convolution and the bound $\norma{v^\nu - \tilde
    v^\nu}_{\L1} \leq \O \int_{\reali}\sum_{i=1}^n \modulo{q_i^-(x)}
  dx$ yield
  \begin{eqnarray*}
    & &
    \mathbf{\Phi} \left( \F_N(s) v^\nu, \F_N(s) \tilde v^\nu \right)
    -
    \mathbf{\Phi} (v^\nu, \tilde v^\nu)
    \\
    &\leq &
    \left( -\frac{c}{2} +\O\delta\right) s
    \int_{\reali} \sum_{i=1}^n \modulo{q_i^-(x)} dx
    \\
    & &
    \qquad
    + 
    \O s
    \norma{\Pi_N Q*v^\nu - \Pi_N Q*\tilde v^n}_{\L1}
    \\
    & \leq &
    \left( -\frac{c}{2} +\O\delta\right) s
    \int_{\reali} \sum_{i=1}^n \modulo{q_i^-(x)} dx
    \\
    & &
    \qquad
    + 
    \O s \norma{Q}_{\L1} \int_{\reali} \sum_{i=1}^n \modulo{q_i^-(x)} dx
    \\
    & \leq &
    \left( -\frac{c}{2} +\O \left( \delta + \norma{Q}_{\L1} \right) \right) 
    s
    \int_{\reali} \sum_{i=1}^n \modulo{q_i^-(x)} dx
    \\
    & \leq &
    -\frac{c}{4} s
    \int_{\reali} \sum_{i=1}^n \modulo{q_i^-(x)} dx
    \\
    & \leq &
    -\frac{c}{4} \, s \, \mathbf{\Phi} (v^\nu, \tilde v^\nu) \,,
  \end{eqnarray*}
  for $\delta + \norma{Q}_{\L1}$ sufficiently small. The proof is
  completed by means of the lower semicontinuity of $\mathbf{\Phi}$:
  \begin{eqnarray*}
    \mathbf{\Phi} \left(\F(s) v, \F(s) \tilde v \right)
    & \leq &
    \liminf_{N\to +\infty} \liminf_{\nu \to +\infty}
    \mathbf{\Phi} \left( \F_N(s) v^\nu, \F_N(s) \tilde v^\nu \right)
    \\
    & \leq &
    \lim_{\nu \to +\infty} 
    \left( 1 - \frac{c}{4} s \right)
    \mathbf{\Phi} \left( v^\nu, \tilde v^\nu \right)
    \\
    & = &
    \left( 1 - \frac{c}{4} s \right)
    \mathbf{\Phi} \left( v, \tilde v \right) .
  \end{eqnarray*}
\end{proofof}

\begin{proofof}{Theorem~\ref{thm:BL}}
  Let $\delta > 0$ be so small that lemmas~\ref{lem:V}
  and~\ref{lem:Stime} hold. First, we show that Theorem~\ref{thm:main}
  can be applied with $\hat \mathcal{D}= \mathcal{D}^M$ and $F$ as
  in~\Ref{eq:LocalFlow}.

  For $u \in \mathcal{D}$, the map $G$ defined in~\Ref{eq:source} is
  $\L1$-bounded, $\L1$-Lipschitz and $\tv\left( G(u) \right)$ is
  uniformly bounded.  The Lipschitz constant of $F$ with respect to
  time can be estimated as follows:
  \begin{eqnarray}
    \nonumber
    & &
    \norma{F(s) u - F(s')u}_{\L1}
    \\
    \nonumber
    & = &
    \norma{S_s u + s G(S_s u) - S_{s'} u - s' G(S_{s'}u)}_{\L1}
    \\
    \nonumber
    & \leq &
    L \, \modulo{s-s'}
    +
    \modulo{s-s'} \, \norma{G(S_s u)}_{\L1}
    +
    s' \, \norma{G(S_s u) - G(S_{s'} u)}_{\L1}
    \\
    \label{eq:LipT}
    & \leq &
    \O \modulo{s-s'} \left( 1+\norma{u}_{\L1} \right)
  \end{eqnarray}
  hence $F$ is Lipschitz in $t$ uniformly in $u\in \mathcal{D}^M$, by
  the boundedness of $\mathcal{D}^M$ in $\L1$. The Lipschitzeanity in
  $u$ is straightforward.

  $F$ satisfies condition~\ref{it:first}, indeed
  \begin{eqnarray*}
    & &
    \norma{F(ks) F(s) u - F \left((k+1)s\right)u}_{\L1}
    \\
    & \leq &
    \norma{
      S_{ks} \left( S_s u + sG(S_su) \right) 
      - 
      S_{ks} S_s u
      -
      s G(S_{ks} S_s u)
    }_{\L1}
    \\
    & &
    +
    ks 
    \norma{
      G\left(S_{ks} \left( S_s u + sG(S_s u) \right) \right)
      - G(S_{ks} S_s u)
    } \,.
  \end{eqnarray*}
  Apply~\cite[Proposition~3.10]{ColomboGuerra} to the first term:
  \begin{eqnarray*}
    & &
    \norma{
      S_{ks} \left( S_s u + sG(S_su) \right) 
      - 
      S_{ks} S_s u
      -
      s G(S_{ks} S_s u)
    }_{\L1}
    \\
    & \leq &
    \O \norma{sG(S_su) - s G(S_{ks} S_s u)}_{\L1}
    +
    \O \, ks \, \tv \left( s G(S_{ks} S_s u) \right)
    \\
    & \leq &
    \O \, ks \, s \,.
  \end{eqnarray*}
  The second term is of the same order by the Lipschitzeanity of $G$
  and of the SRS. Therefore, condition~\ref{it:first} is proved with
  $\omega(s) = s$.

  $F$ satisfies the stability condition~\ref{it:second}. Indeed, by
  Lemma~\ref{lem:Stime}
  \begin{displaymath}
    \mathbf{\Phi} \left( F^\epsilon(t)u, F^\epsilon(t)\tilde u \right)
    \leq
    e^{-(c/4)t} \, \mathbf{\Phi}(u, \tilde u) \,.
  \end{displaymath}
  By~3 in Proposition~\ref{prop:Functionals},
  \begin{equation}
    \label{eq:Lipschitz}
    \norma{F^\epsilon(t)u - F^\epsilon(t)\tilde u}_{\L1}
    \leq
    \mathcal{L} \, e^{-(c/4)t} \,
    \norma{u-\tilde u}_{\L1}
  \end{equation}
  with $\mathcal{L}$ independent from $\epsilon$ and $M$.

  Applying Theorem~\ref{thm:main}, we obtain for all $u \in
  \mathcal{D}^M$ the strong $\L1$ convergence $F^\epsilon(t) u \to P_t
  u$, $P$ being the unique $\L1$-Lipschitz semigroup
  satisfying~\Ref{eq:tangent}. Moreover, passing to the limit
  $\epsilon \to 0+$ in~\Ref{eq:Lipschitz}, we obtain the first
  estimate in~\Ref{eq:error}, with $\mathcal{L}$ independent from
  $M$. Therefore, letting $M \to + \infty$, we may uniquely extend $P$
  to all $\mathcal{D}_\delta$, keeping the validity of the first
  estimate in~\Ref{eq:error}.

  Fix $u \in \mathcal{D}_\delta$ and let $M_u = \mathbf{\Phi}(u,0)$,
  so that $u \in \mathcal{D}^{M_u}$. By~\Ref{eq:LipT}
  and~\cite[Lemma~2.3]{ColomboGuerra3},
  \begin{eqnarray*}
    \norma{F^\epsilon (t)u - F^\epsilon (s)u}_{\L1} 
    & \leq &
    \O \, ( 1 + M_u ) \, \modulo{t-s}
    \\
    & \leq &
    \O \, \left( 1 + \norma{u}_{\L1} \right) \, \modulo{t-s}
  \end{eqnarray*}
  and passing to the limit $\epsilon \to 0$ we obtain the second
  estimate in~\Ref{eq:error}.

  Condition~(4) is a direct consequence of~\Ref{eq:tangent} and it
  ensures that the orbits of $P$ are weak entropy solutions,
  see~\cite[Corollary~3.13]{ColomboGuerra}.

  Finally, \Ref{eq:Uffa} is obtained as
  in~\cite[formula~(1.7)]{ColomboGuerra}, see
  also~\cite{BianchiniColombo}.
\end{proofof}

{\small{

    \bibliographystyle{abbrv}

  }}

\end{document}